\documentclass[11pt,english]{article}
\usepackage[T1]{fontenc}
\usepackage[latin9]{inputenc}
\usepackage{geometry}
\usepackage{amsmath}
\usepackage{amssymb}
\usepackage{graphicx}
\usepackage{setspace}
\makeatletter

\usepackage{tikz-cd}

\makeatother

\usepackage{babel}
\begin{document}

\title{On the $v_{1}$-periodicity of the Moore space}

\author{Lyuboslav Panchev}
\maketitle
\begin{abstract}
We present progress in trying to verify a long-standing conjecture
by Mark Mahowald on the $v_{1}$-periodic component of the classical
Adams spectral sequence for a Moore space $M$. The approach we follow
was proposed by John Palmieri in his work on the stable category of
$A$-comodules. We improve on Palmieri's work by working with the
endomorphism ring of $M$ - $End(M)$ thus resolving some of the initial
difficulties of his approach and formulating a conjecture of our own
that would lead to Mahowald's formulation. 
\end{abstract}

\section{Introduction}

Stable homotopy groups of spheres (or generally any finite complex)
have long been a subject of study in algebraic topology. Due to their
immense complexity, one might try to understand them from the lens
of the Adams spectral sequence. However, this task as well seems outside
of the scope of what we can fully understand. Nevertheless, the chromatic
point of view allows us to ``break'' the spectral sequence into
chromatic pieces. This paper represents an attempt to understand one
of those pieces - namely, the $v_{1}$-periodic component of the Adams
spectral sequence of a Moore space i.e. $v_{1}^{-1}E_{2}(M;H)$. A
conjecture, which seems quite likely given the available data, was
formulated in \cite{key-1} by Mark Mahowald almost 50 years ago.
Still, of course, no amount of data can be a substitute for a rigorous
proof.

Palmieri proposed an approach to this conjecture in \cite{key-2}.
He built a generalized Adams spectral sequence in the stable category
of comodules over the dual Steenrod algebra $A$. This spectral sequence
converges to $v_{1}^{-1}E_{2}(M;H)$ and computations seem promising
due to the simplicity of $E_{2}=\mathbb{F}_{2}[v_{1}^{\pm1},h_{11},h_{21},\cdots,h_{n1},\cdots]$
and the fact that $E_{3}=E_{2}$ as for degree reasons nontrivial
differentials can only occur at odd pages. It is important to note
that since $M$ is not a ring spectrum, $E_{r}$ is not an algebra
and $d_{r}$ is not a derivation and so what we really mean by the
above equality is that $E_{2}$ is a $\mathbb{F}_{2}$-vector space
with basis the monomials in $\mathbb{F}_{2}[v_{1}^{\pm1},h_{11},h_{21},\cdots,h_{n1},\cdots]$.
Palmieri then conjectured what the values of $d_{3}(h_{n1})$ are
and proposed one should be able to extend them in some way to the
entire $E_{3}$. Moreover he conjectured that the spectral sequence
collapses at $E_{4}$ and claimed this would imply Mahowald's conjecture.
Note it is not immediately obvious how Palmieri's formulation relates
to Mahowald's and it is something we address in more detail at a later
section of the paper.

Thus our problem is three-fold: how does one compute $d_{3}(h_{n1})$,
how does one extend it to the rest of $E_{3}$ and why are there no
higher degree differentials. We solely address the first two questions,
fully answering the second one. We do this by working with the endomorphism
ring spectrum of $M$ - $End(M)$. It is the $4$ cell complex $M\wedge DM$.
The advantage of $End(M)$ is that its spectral sequence is multiplicative
and so $d_{3}$ is a derivation. At the same time the action $End(M)\wedge M\to M$
makes $E_{r}(M)$ into a module over $E_{r}(End(M))$. We will also
show Palmieri's originally conjectured values for $d_{3}(h_{n1})$
can't be true and so we propose a revised conjecture of what those
values are. We verify that conjecture modulo knowing that the elements
$v_{1}^{m}h_{n1}$ don't survive to $E_{4}$ for $n\ge3$, $m\in\mathbb{Z}$.

In section $2$ we provide the necessary background about $Stable(A)$
- the stable category of comodules over the Steenrod algebra $A$,
and explicitly write Palmieri's original conjecture and our revised
version of it. In section 3 we work out the corresponding spectral
sequence for $End(M)$ and its action on the the one for $M$. Section
$4$ consists of the meat of the paper as we proceed to prove our
main results stated above. We conclude with section $5$ where we
introduce the original conjecture by Mahowlad and show explicitly
how it follows from our revised conjecture.

The author would also like to thank his advisor Haynes Miller for
the unending support and multitude of fruitful discussions and suggestions,
including the crucial idea of working with $End(M)$.

\section{The category $Stable(A)$}

In this chapter we give a brief description of $Stable(A)$ and any
related results of immediate use to us. For more detail the reader
is directed to Palmieri's book \cite{key-2}. 

Objects in $Stable(A)$ are unbounded cochain complexes of (left)
$A$-comodules. We will identify a comodule $L$ with its injective
resolution over $A$. For two such objects $L,N$ the set of morphisms
is $[L,N]_{s,t}=Ext_{A}^{s,t}(L,N)$. Then $L_{s,t}=\pi_{s,t}(L)=Ext_{A}^{s,t}(\mathbb{F}_{2},L)$.
For the sake of clarity we observe $L$ itself is bigraded and one
should make a distinction between the elements of degree $(s,t)$
in $L$ and $L_{s,t}$. Note also the sphere spectrum $S\in Stable(A)$
is the injective resolution of $\mathbb{F}_{2}$, which is in line
with our notation of $\pi_{s,t}(L)=[S,L]_{s,t}$ above. $Stable(A)$
is now a triangulated category and for a ring spectrum $X\in Stable(A)$
we can build a generalaized Adams spectral sequence in the usual way.
Then assuming certain conditions hold we can identify $E_{2}(L;X)=Ext_{X_{**}X}(X_{**},X_{**}L)$
and further conditions would guarantee convergence to $\pi_{**}L$. 

We are interested in the case where the spectrum $Q_{1}$ plays the
role of $X$. To define $Q_{1}$, we first define $q_{1}$ to be the
injective resolution of $A\Box_{\mathbb{F}_{2}(\xi_{2})/(\xi_{2}^{2})}\mathbb{F}_{2}$.
$Q_{1}$ is now obtained from $q_{1}$ after working out how to extend
the $q_{1}$-resolution into the negative dimensions. Then one can
check $q_{1**}=\mathbb{F}_{2}[v_{1}]$, $Q_{1**}=\mathbb{F}_{2}[v_{1}^{\pm1}]$
\cite[p.44]{key-2} and $Q_{1**}Q_{1}=\mathbb{F}_{2}[v_{1}^{\pm1},\xi_{1},\xi_{2}^{2},\cdots\xi_{n}^{2},\cdots]/(\xi_{1}^{4},\xi_{2}^{4},\cdots)$
\cite[p.101]{key-2}.

The trigraded spectral sequence of interest is

\[
E_{2}(M;Q_{1})=Ext_{Q_{1**}Q_{1}}(Q_{1**},Q_{1**}(M))=\mathbb{F}_{2}[v_{1}^{\pm1},h_{11},h_{21},\cdots,h_{n1},\cdots]
\]
and it converges to $v_{1}^{-1}E_{2}(M;H)$ \cite[p.81, 101]{key-2}.
Note the abuse of notation above as what we really mean by $E_{2}(M;Q_{1})$
is $E_{2}(L;Q_{1})$ where $L$ is an injective resolution for $H_{*}(M)$.
Elsewhere $M$ will always refer to the topological Moore spectrum.
For degree reasons the only potential non-zero differentials in $E_{r}(M;Q_{1})$
happen at odd pages, so $E_{2}=E_{3}$. Palmieri then conjectured
the following differentials:

\[
\begin{array}{cc}
d_{3}(v_{1}^{2})=h_{11}^{3}\\
d_{3}(h_{n1})=v_{1}^{-2}h_{11}h_{21}h_{n-1,1}^{2} & \text{for }n\ge3
\end{array}
\]

As we will see later, the conjecture in its current form is incorrect,
so we make the following revised conjecture:

\[
\begin{array}{cc}
d_{3}(v_{1}^{2})=h_{11}^{3}\\
d_{3}(h_{n1})=v_{1}^{-2}h_{11}^{3}h_{n1}+v_{1}^{-2}h_{11}h_{21}h_{n-1,1}^{2} & \text{for }n\ge3
\end{array}
\]
Though this isn't enough to fully determine $d_{3}$, Palmieri goes
on to propose that $d_{3}$ ``looks'' as though as $E_{2}(M;Q_{1})$
is an algebra. One reason for this proposal that he notes is we can
also compute the $E_{2}$ page of the corresponding spectral sequence
for the sphere 
\[
E_{2}(S;Q_{1})=Ext_{Q_{1**}Q_{1}}(Q_{1**},Q_{1**})=\mathbb{F}_{2}[v_{1}^{\pm1},h_{10},h_{11},h_{21},\cdots,h_{n1},\cdots]
\]
and use the map $S\to M$ to induce a surjection $E_{2}(S;Q_{1})\to E_{2}(M;Q_{1})$
with $h_{n1}\to h_{n1}$, $h_{10}\to0$ and $v_{1}\to v_{1}$. Then
the identity map $S\wedge M\to M$ turns $E_{2}(M;Q_{1})$ into a
cyclic module over $E_{2}(S;Q_{1})$. Now identifying $E_{2}(M;Q_{1})$
with $\mathbb{F}_{2}[v_{1}^{\pm1},h_{11},h_{21},\cdots,h_{n1},\cdots]$
becomes justified as both coincide as $E_{2}(S;Q_{1})$-modules: 
\[
E_{2}(M;Q_{1})\cong E_{2}(S;Q_{1})/(h_{10})=\mathbb{F}_{2}[v_{1}^{\pm1},h_{11},h_{21},\cdots,h_{n1},\cdots]
\]

Then information about differentials in $E_{r}(S;Q_{1})$ could directly
produce differentials in $E_{r}(M;Q_{1})$ and since $S$ is a ring
spectrum, $E_{r}(S;Q_{1})$ is a spectral sequence of algebras, so
the differentials in $E_{r}(S;Q_{1})$ are derivations. The problem
is differentials in $E_{2}(S;Q_{1})$ are difficult to compute and
so we don't know what $E_{3}(S;Q_{1})$ looks like. This is where
$End(M)$ enters the picture - it is a ring spectrum that acts on
$M$ just as $S$ does, but differentials in $E_{2}(End(M);Q_{1})$
are much more manageable to compute.

\section{The $Q_{1}$ $E_{2}$ term for $End(M)$}

We begin by computing $H_{*}(End(M))$ as a comodule over $A$. Let
$x_{0}$ and $x_{1}$ denote the two cells of $M$ and $y_{-1}$ and
$y_{0}$ denote the two cells of $DM=\Sigma^{-1}M$. Then $End(M)=M\wedge DM$
has four cells of the form $x_{i}y_{j}$ with $|x_{i}y_{j}|=i+j$.
As $DM$ is the dual of $M$ we have maps $\eta:S\to M\wedge DM$
and $\epsilon:DM\wedge M\to S$ that specify the ring structure of
$End(M)$. More precisely, $\eta$ is the unit, while multiplication
is given by\[\begin{tikzcd} 
	M\wedge DM\wedge M\wedge DM \arrow[r,"1\wedge \epsilon\wedge 1"] & M\wedge DM 
\end{tikzcd}\]

\noindent and the action of $End(M)$ on $M$ is then given by the
map $1\wedge\epsilon:M\wedge DM\wedge M\to M$. If $\iota\in H_{*}(S)$
is the generator, then $\eta_{*}(\iota)=x_{1}y_{-1}+x_{0}y_{0}$ and
$\epsilon_{*}(y_{1}x_{-1})=\epsilon_{*}(y_{0}x_{0})=\iota$. This
allows us to compute the multiplicative structure of $H_{*}(End(M))$
\[
(x_{i}y_{j})(x_{k}y_{l})=\begin{cases}
\begin{array}{cc}
x_{i}y_{l} & \text{if }j+k=0\\
0 & \text{otherwise}
\end{array}\end{cases}
\]

Setting $\alpha=x_{0}y_{-1}$ and $\gamma=x_{1}y_{0}$ we get that
$H_{*}(End(M))=\mathbb{F}_{2}[\alpha,\gamma]/(\alpha^{2},\gamma^{2},\alpha\gamma+\gamma\alpha+1)$.
Note this is a 4-dimensional non-commutative $\mathbb{F}_{2}$-algebra
with basis $\langle1,\alpha,\gamma,\alpha\gamma\rangle$ where $|\alpha|=-1$
and $|\gamma|=1$. To understand the coaction of $A$ we just need
to understand the coaction on $\alpha$ and $\gamma$. Since $\psi(x_{0})=1\otimes x_{0}$
and $\psi(x_{1})=1\otimes x_{1}+\xi_{1}\otimes x_{0}$ we conclude
that 
\[
\psi(\alpha)=\psi(x_{0}y_{-1})=\psi(x_{0})\psi(y_{-1})=(1\otimes x_{0})(1\otimes y_{-1})=1\otimes x_{0}y_{-1}=1\otimes\alpha
\]
 and 
\begin{eqnarray*}
\psi(\gamma) & = & \psi(x_{1}y_{0})=\psi(x_{1})\psi(y_{0})=1\otimes x_{1}y_{0}+\xi_{1}\otimes(x_{1}y_{-1}+x_{0}y_{0})+\xi_{1}^{2}\otimes x_{0}y_{-1}\\
 & = & 1\otimes\gamma+\xi_{1}\otimes1+\xi_{1}^{2}\otimes\alpha
\end{eqnarray*}

Recall we are interested in computing $d_{3}$ in $E_{2}(M;Q_{1})$.
Since $M$ lacks multiplicative structure, we will work with $End(M)$
and try to understand $E_{r}(End(M);Q_{1})$. We proceed with a direct
computation 
\begin{eqnarray*}
E_{2}(End(M);Q_{1}) & = & Ext_{(Q_{1})_{**}Q_{1}}((Q_{1})_{**},(Q_{1})_{**}(End(M)))\\
 & = & \mathbb{F}_{2}[v_{1}^{\pm1}]\otimes Ext_{\mathbb{F}_{2}[\xi_{1},\xi_{2}^{2},\cdots]/(\xi_{i}^{4})}(\mathbb{F}_{2},\mathbb{F}_{2}\langle1,\alpha,\gamma,\alpha\gamma\rangle)\\
 & = & \mathbb{F}_{2}[v_{1}^{\pm1}]\otimes\mathbb{F}_{2}[h_{21},h_{31},...]\otimes Ext_{\mathbb{F}_{2}[\xi_{1}]/(\xi_{1}^{4})}(\mathbb{F}_{2},\mathbb{F}_{2}\langle1,\alpha,\gamma,\alpha\gamma\rangle)
\end{eqnarray*}

Here we used that the coaction of $\xi_{i}^{2}$ on $\mathbb{F}_{2}\langle1,\alpha,\gamma,\alpha\gamma\rangle$
is trivial for $i\ge2$. The conormal extension $\mathbb{F}_{2}(\xi_{1}^{2})/(\xi_{1}^{4})\to\mathbb{F}_{2}(\xi_{1})/(\xi_{1}^{4})\to\mathbb{F}_{2}(\xi_{1})/(\xi_{1}^{2})$
produces a Cartan-Eilenberg spectral sequence that collapses since
$H_{*}(End(M))=\mathbb{F}_{2}\langle1,\alpha,\gamma,\alpha\gamma\rangle$
is cofree over $\mathbb{F}_{2}(\xi_{1})/(\xi_{1}^{2})$. Thus, we
get

\begin{eqnarray*}
Ext_{\mathbb{F}_{2}[\xi_{1}]/(\xi_{1}^{4})}(\mathbb{F}_{2},\mathbb{F}_{2}\langle1,\alpha,\gamma,\alpha\gamma\rangle) & = & Ext_{\mathbb{F}_{2}[\xi_{1}^{2}]/(\xi_{1}^{4})}(\mathbb{F}_{2},Ext_{\mathbb{F}_{2}[\xi_{1}]/(\xi_{1}^{2})}(\mathbb{F}_{2},\mathbb{F}_{2}\langle1,\alpha,\gamma,\alpha\gamma\rangle))\\
 & = & Ext_{\mathbb{F}_{2}[\xi_{1}^{2}]/(\xi_{1}^{4})}(\mathbb{F}_{2},\mathbb{F}_{2}\langle1,\alpha\rangle)\\
\end{eqnarray*}

We conclude that $Ext_{\mathbb{F}_{2}[\xi_{1}]/(\xi_{1}^{4})}(\mathbb{F}_{2},\mathbb{F}_{2}\langle1,\alpha,\gamma,\alpha\gamma\rangle)=\mathbb{F}_{2}\langle1,\alpha\rangle\otimes\mathbb{F}_{2}[h_{11}]$
and so 
\[
E_{2}(End(M);Q_{1})=\mathbb{F}_{2}[v_{1}^{\pm1},\alpha,h_{11},h_{21},h_{31},...]/(\alpha^{2})
\]
which (expectedly so) is two copies of $E_{2}(M;Q_{1})$. The degrees
of the generators are given by $|v_{1}|=(0,2,1),\,|\alpha|=(0,-1,0),\,|h_{n1}|=(1,2^{n+1}-2,0)$.
It is worth noting that even though $H_{*}(End(M))$ is not commutative,
the spectral sequence above ends up with a commutative multiplicative
structure.

\subsection{$E_{2}(M;Q_{1})$ as a differential module over $E_{2}(End(M);Q_{1})$}

The action of $End(M)$ on $M$ extends to an action $E_{r}(End(M);Q_{1})\otimes E_{r}(M;Q_{1})\to E_{r}(M;Q_{1})$
and so $E_{r}(M;Q_{1})$ is a differential module over $E_{r}(End(M);Q_{1})$.
The commutative diagram

\[\begin{tikzcd} 
	M\wedge DM \wedge M \arrow[r,"1\wedge \epsilon"] & M \\
	S \wedge M \arrow[u,"\eta\wedge 1"] \arrow[ru, "\cong"]&
\end{tikzcd}\]

\noindent implies the action of $E_{2}(S;Q_{1})$ on $E_{2}(M;Q_{1})$
factors through the action of $E_{2}(End(M);Q_{1})$ via the algebra
map $\eta_{*}:E_{r}(S,Q_{1})\to E_{r}(End(M);Q_{1})$, which is just
\[
\eta_{*}:\,\mathbb{F}_{2}[v_{1}^{\pm1}]\otimes\mathbb{F}_{2}[h_{10},h_{11},h_{21},h_{31},...]\to\mathbb{F}_{2}[v_{1}^{\pm1}]\otimes\mathbb{F}_{2}[h_{11},h_{21},h_{31},...]\otimes\mathbb{F}_{2}\langle1,\alpha\rangle
\]

\noindent with $\eta_{*}(v_{1})=v_{1}$ and $\eta_{*}(h_{n1})=h_{n1}$.
Furthermore we claim $\eta_{*}(h_{10})=\alpha h_{11}$. Indeed, since
$\psi(\gamma)=1\otimes\gamma+\xi_{1}\otimes1+\xi_{1}^{2}\otimes\alpha$
it follows that $\xi_{1}\otimes1+\xi_{1}^{2}\otimes\alpha$ vanishes
in the homology of the cobar complex of $End(M)$ and so $\alpha h_{11}=\xi_{1}^{2}|\alpha=\xi_{1}|1$,
which is the cobar representative of $h_{10}$ in $E_{2}(S;Q_{1})$.

Hence $E_{2}(M;Q_{1})$ is a cyclic module over $E_{2}(End(M);Q_{1})$.
Furthermore,  we have an isomorphism of $E_{2}(End(M);Q_{1})$-modules:

\[
E_{2}(M;Q_{1})\cong E_{2}(End(M);Q_{1})/(\alpha)=\mathbb{F}_{2}[v_{1}^{\pm1},h_{11},h_{21},\cdots,h_{n1},\cdots]
\]

Before we move on to the next section we note that all of the elements
$h_{11},v_{1},h_{21}v_{1},h_{21}v_{1}^{2}$ survive to $E_{\infty}(M;Q_{1})$
as shown by the diagram of $E_{2}(M;H)$ below. Observe this doesn't
guarantee the same is true in $E_{r}(End(M);Q_{1}),$ but we will
still be able to extract some of the information back to $E_{r}(End(M);Q_{1})$
using the action above.

~

\bigskip{}

\qquad{} \quad{}\includegraphics{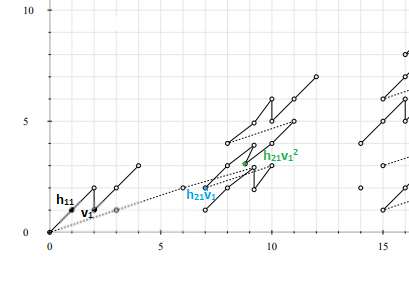}

\section{Calculating $d_{2}$ and $d_{3}$ of $E_{2}(End(M);Q_{1})$}

We begin by calculating $d_{2}$ and $d_{3}$ on the low-degree elements
in $E_{r}(End(M);Q_{1})$ and then proceed to formulating a conjecture
for $d_{2}$ and $d_{3}$ on the remaining elements. 

$\textbf{Theorem 1:}$ The elements $\alpha,h_{11},v_{1}\alpha,v_{1}h_{21}$
survive to $E_{4}(End(M);Q_{1})$. Furthermore,

$d_{2}(v_{1})=\alpha h_{11}^{2}$

$d_{3}(v_{1}^{2})=h_{11}^{3}$

$\emph{Proof of Theorem 1:}$

Since we will need to distinguish between differentials in $E_{r}(End(M);Q_{1})$
and $E_{r}(M;Q_{1})$, we will denote them by $d_{r}$ and $d_{r}^{M}$
respectively. 

In $E_{r}(M;Q_{1})$, $h_{11}^{3}$ must be a coboundary at some point
and for degree reasons $d_{3}^{M}(v_{1}^{2})=h_{11}^{3}$. Indeed,
if $d_{r}(x)=h_{11}^{3}$ for some $r\ge3$ and $x\in E_{r}(M;Q_{1})$
then since $|h_{11}^{3}|=(3,6,0)$ and $d_{r}^{M}$ changes degrees
by $(r,r-1,1-r)$ we conclude that $|x|=(3-r,7-r,r-1)$. Recall $|v_{1}|=(0,2,1),\,|\alpha|=(0,-1,0),\,|h_{n1}|=(1,2^{n+1}-2,0)$.
Then $3-r\ge0$, so $r=3$ and $|x|=(0,4,2)$. The only option now
is $x=v_{1}^{2}$. Note if $v_{1}$ was to survive to $E_{3}(End(M);Q_{1})$
then $d_{3}(v_{1}^{2})=0$, which would force $d_{3}^{M}(v_{1}^{2})=0$.
Hence $d_{2}(v_{1})\neq0$ and so for degree reasons $d_{2}(v_{1})=\alpha h_{11}^{2}$.
Given the action of $E_{2}(End(M);Q_{1})$ we must also have $d_{2}(v_{1})=\alpha h_{11}^{2}$.
Either of those differentials could be also seen since $d_{2}(v_{1})=h_{10}h_{11}$
in $E_{2}(S;Q_{1})$ which follows from the same differential in the
Cartan-Eilenberg spectral sequence computing $H^{*}(A(1))$.

Next we claim $d_{2}(h_{21})\neq0$. Indeed, assume that $d_{2}(h_{21})=0$.
Then $d_{2}(v_{1}^{2}h_{21})=0$ and since $v_{1}^{2}h_{21}$ survives
in $E_{r}(M;Q_{1})$ it must be that $d_{3}(v_{1}^{2}h_{21})=0$ in
$E_{3}(End(M);Q_{1})$. By multiplicativity we conclude $d_{3}(h_{21})=v_{1}^{-2}h_{11}^{3}h_{21}$.
But now considering the action $E_{3}(End(M);Q_{1})\otimes E_{3}(M;Q_{1})\to E_{3}(M;Q_{1})$
we have 
\[
d_{3}^{M}(h_{21}\cdot v_{1})=d_{3}(h_{21})\cdot v_{1}+h_{21}\cdot d_{3}^{M}(v_{1})=v_{1}^{-1}h_{11}^{3}h_{21}\neq0
\]
which can't happen since $h_{21}v_{1}$ survives in $E_{r}(M;Q_{1})$.
Note we have to consider the action since $h_{21}v_{1}$ would not
be present in $E_{3}(End(M);Q_{1})$. Hence our assumption was wrong
and $d_{2}(h_{21})\neq0$, which by degree reasons means $d_{2}(h_{21})=v_{1}^{-1}\alpha h_{11}^{2}h_{21}$.

Finally both $h_{11}$ and $v_{1}h_{21}$ survive $d_{3}^{M}$ in
$E_{3}(M;Q_{1})$, so they must also survive $d_{3}$ in $E_{3}(End(M);Q_{1})$
i.e. $d_{3}(h_{11})=d_{3}(v_{1}h_{21})=0$. At the same time, for
degree reasons $d_{r}(\alpha)=d_{r}(\alpha v_{1})=0$ for $r=2,3$
and neither elements can be a coboundary, which means both $\alpha$
and $\alpha v_{1}$ are present in $E_{4}(End(M);Q_{1})$.

\hfill{}$\Box$

Given the theorem above, in order to compute $d_{2}$ completely we
just need to know the values on the remaining generators i.e. $d_{2}(h_{n1})$
for $n\ge3$. Thus we make the following conjecture: 
\[
\textbf{(Main) Conjecture part 1: }d_{2}(h_{n1})=v_{1}^{-1}\alpha h_{11}^{2}h_{n1}\text{ for }n\ge3\,\,\,\,\,\,\,\,\,\,\,\,\,\,\,\,\,\,\,\,\,\,
\]
. 

\medskip{}

\noindent Observe then $x_{n}=v_{1}h_{n+1,1}$ is a cycle, and that

\[
E_{2}(End(M);Q_{1})=\mathbb{F}_{2}[x_{1},x_{2},...]\otimes\mathbb{F}_{2}[v_{1}^{\pm1},h_{11},\alpha]/(\alpha^{2})
\]

\noindent where the first factor has zero differential and the second
factor has only $d_{2}v_{1}=\alpha h_{11}^{2}$. The homology is thus 

\[
E_{3}(End(M);Q_{1})=\mathbb{F}_{2}[x_{1},x_{2},...]\otimes\mathbb{F}_{2}[v_{1}^{\pm2},h_{11},\alpha,\alpha']/(\alpha^{2},\alpha h_{11}^{2},\alpha\alpha',\alpha'^{2})
\]

\noindent where $\alpha'$ is the class of $v_{1}\alpha$. Again Theorem
1 tells us $d_{3}(x_{1})=d_{3}(\alpha)=d_{3}(\alpha')=0$ and $d_{3}(v_{1}^{2})=h_{11}^{3}$
and so in order to compute $d_{3}$ completely we just need to know
the values on the remaining generators i.e. $d_{3}(x_{n})$ for $n\ge2$.
Thus we further conjecture:

\[
\textbf{(Main) Conjecture part 2: }d_{3}(x_{n})=v_{1}^{-4}h_{11}x_{1}x_{n-1}^{2}\text{ for }n\ge2\,\,\,\,\,\,\,\,\,\,\,\,\,\,\,\,\,\,\,\,\,
\]

\smallskip{}

\noindent We can prove this conjecture modulo the following assumption:

$\textbf{(Smaller) conjecture:}$ $v_{1}^{m}x_{n}$ does not survive
to $E_{4}(End(M);Q_{1})$ for $n,m\in\mathbb{Z}$, $n\ge2$.\medskip{}

$\textbf{Theorem 2:}$ The smaller conjecture above implies the main
one.

Before proving the Theorem observe the converse statement that the
main conjecture implies the smaller one also holds. In fact, the main
conjecture even specifies what $d_{r}(v_{1}^{m}x_{n})$ is, which
is what justifies the naming convention of the two conjectures. Thus,
the Theorem can be reformulated by saying that the smaller and main
conjectures above are equivalent.

$\emph{Proof of Theorem 2:}$

For $n\ge3$ $d_{2}(h_{n1})$ is a linear combination of $v_{1}^{-1}\alpha h_{11}^{2}h_{n1}$
and $v_{1}^{-1}\alpha h_{21}h_{n-1,1}^{2}$ for degree reasons, but
the later is not in the image of $E_{2}(S;Q)$. Hence $d_{2}(h_{n1})=v_{1}^{-1}\alpha h_{11}^{2}h_{n1}$
or $0$. Assume that for some $n\ge3$ $d_{2}(h_{n1})=0$. For degree
reasons, $d_{3}(h_{n1})$ is a linear combination of $v_{1}^{-2}h_{11}^{3}h_{n1}$
and $v_{1}^{-2}h_{11}h_{21}h_{n-1,1}^{2}$, but $v_{1}^{-2}h_{11}h_{21}h_{n-1,1}^{2}$
doesn't survive to $E_{3}(End(M);Q_{1})$ since 
\[
d_{2}(v_{1}^{-2}h_{11}h_{21}h_{n-1,1}^{2})=d_{2}(h_{21})v_{1}^{-2}h_{11}h_{n-1,1}^{2}=v_{1}^{-3}\alpha h_{11}^{3}h_{21}h_{n-1,1}^{2}
\]
By our smaller conjecture, $d_{3}(h_{n1})\neq0$ and so $d_{3}(h_{n1})=v_{1}^{-2}h_{11}^{3}h_{n1}$.
Then 
\[
d_{3}(v_{1}^{2}h_{n1})=d_{3}(v_{1}^{2})h_{n1}+v_{1}^{2}d_{3}(h_{n1})=h_{11}^{3}h_{n1}+h_{11}^{3}h_{n1}=0
\]
which again contradicts the (smaller) conjecture. We conclude $d_{2}(h_{n1})=v_{1}^{-1}\alpha h_{11}^{2}h_{n1}$
for all $n\ge2$, which is also equivalent to $d_{2}(v_{1}h_{n1})=0$
for all $n\ge2$. Hence the elements $x_{n}=v_{1}h_{n+1,1}$ survive,
which justifies their presence in $E_{3}$. This completes the $d_{2}$
calculation in $E_{2}(End(M);Q_{1})$.

Next for $n\ge2$ $d_{3}(x_{n})$ is a linear combination of $v_{1}^{-4}h_{11}x_{1}x_{n-1}^{2}$
and $v_{1}^{-2}h_{11}^{3}x_{n}$, which leaves us with 4 possibilities.
$d_{3}(x_{n})=v_{1}^{-2}h_{11}^{3}x_{n}$ would imply $d_{3}(v_{1}^{2}x_{n})=0$
and so $d_{3}(x_{n})=0$ or $v_{1}^{-2}h_{11}^{3}x_{n}$ are both
ruled out as possibilities due to the (smaller) conjecture. Then either
$d_{3}(x_{n})=v_{1}^{-4}h_{11}x_{1}x_{n-1}^{2}$ or $d_{3}(x_{n})=v_{1}^{-4}h_{11}x_{1}x_{n-1}^{2}+v_{1}^{-2}h_{11}^{3}x_{n}$.
However, the latter case would imply 
\[
d_{3}(v_{1}^{2}x_{n})=d_{3}(v_{1}^{2})x_{n}+v_{1}^{2}d_{3}(x_{n})=h_{11}^{3}x_{n}+h_{11}^{3}x_{n}+v_{1}^{-2}h_{11}x_{1}x_{n-1}^{2}=v_{1}^{-2}h_{11}x_{1}x_{n-1}^{2}
\]

\noindent and so
\[
0=d_{3}^{2}(v_{1}^{2}x_{n})=d_{3}(v_{1}^{-2}h_{11}x_{1}x_{n-1}^{2})=d_{3}(v_{1}^{-2})h_{11}x_{1}x_{n-1}^{2}=v_{1}^{-4}h_{11}^{4}x_{1}x_{n-1}^{2}
\]
which is false as $v_{1}^{-4}h_{11}^{4}x_{1}x_{n-1}^{2}$ is present
in $E_{3}(End(M);Q_{1})$. We conclude $d_{3}(x_{n})=v_{1}^{-4}h_{11}x_{1}x_{n-1}^{2}$
for $n\ge2$ as desired. 

\hfill{}$\Box$

It is worth mentioning that Palmieri's original conjecture would imply
that $d_{3}^{M}(v_{1}^{m}h_{n1})\neq0$ for $n\ge3$, which would
guarantee the (smaller) conjecture. However, the smaller conjecture
itself is enough to arrive at a different answer than what Palmieri
suggested. This proves his original formulation is incorrect, but
as we will see in the next section it is close to what we arrive at
based on the (smaller) conjecture.

\subsection{Completing the calculation of $d_{3}$ in $E_{3}(M;Q_{1})$}

Now that we have learnt a fair bit about the structure of $E_{r}(End(M);Q_{1})$
we will see how the information about its differentials can translate
to information about the differentials in $E_{r}(M;Q_{1})$. Recall
for degree reasons $E_{2}(M;Q_{1})=E_{3}(M;Q_{1})$. Observe $E_{3}(M;Q_{1})$
is now generated by $\{1,v_{1}\}$ as a $E_{3}(End(M);Q_{1})$-module.
Since $v_{1}$ survives to $E_{\infty}(M;Q_{1})$ we get $d_{3}^{M}(v_{1})=d_{3}^{M}(1)=0$
and so $d_{3}$ now completely determines $d_{3}^{M}$. 

For example, to compute $d_{3}^{M}(h_{n1})$ for $n\ge3$ note that
$h_{n1}=v_{1}^{-2}x_{n-1}\cdot v_{1}$ and so we get

\[
d_{3}^{M}(h_{n1})=d_{3}(v_{1}^{-2}x_{n-1})\cdot v_{1}=v_{1}^{-2}h_{11}^{3}h_{n1}+v_{1}^{-2}h_{11}h_{21}h_{n-1,1}^{2}
\]
We conclude that assuming the (smaller) conjecture holds, the differentials
in $E_{3}(M;Q_{1})$ are 

$d_{3}^{M}(v_{1}^{2})=h_{11}^{3}$

$d_{3}^{M}(h_{21})=v_{1}^{-2}h_{21}h_{11}^{3}$

$d_{3}^{M}(h_{n1})=v_{1}^{-2}h_{11}^{3}h_{n1}+v_{1}^{-2}h_{11}h_{21}h_{n-1,1}^{2}$
for $n\ge3$

\noindent which is what we conjectured in Section 2.

\section{Relation between Palmieri's and Mahowald's notations}

In this section we will see how the conjectured differentials for
$E_{3}(M;Q_{1})$ imply Mahowald's conjecture assuming there are no
higher degree differentials. We begin by stating Mahowald's conjecture
explicitly following the original description in \cite{key-1}. Let
$P=\mathbb{F}_{2}[x_{1},x_{2},\cdots]$ be a polynomial algebra, which
is bigraded with $|x_{i}|=(2,2^{i+2}+1)$. Set a derivation $d$ on
$P$ by $d(x_{i})=x_{1}x_{i-1}^{2}$ for $i>1$. Let $H(d)$ be the
resulting homology and $B(d)$ the image of $d$. Then assuming $a$
and $b$ run through an $\mathbb{F}_{2}$-basis for $H(d)$ and $B(d)$
Mahowald conjectured that

\begin{eqnarray*}
\text{ }v_{1}^{-1}Ext_{A}^{s,t}(\mathbb{F}_{2},H_{*}(M)) & = & \bigoplus_{a\in H(d)}\Sigma^{|a|}v_{1}^{-1}Ext_{A}^{s,t}(\mathbb{F}_{2},H_{*}(bo\wedge M))\\
 & \oplus & \bigoplus_{b\in B(d)}\Sigma^{|b|}v_{1}^{-1}Ext_{A}^{s,t}(\mathbb{F}_{2},H_{*}(bu\wedge M))
\end{eqnarray*}

Here $bo$ and $bu$ are connective real and complex $K$-theory respectively
and we have explicit computations: 
\[
v_{1}^{-1}Ext_{A}^{s,t}(\mathbb{F}_{2},H_{*}(bo\wedge M))=\mathbb{F}_{2}[v_{1}^{\pm4}]\otimes\mathbb{F}_{2}(h_{11},v_{1})/(h_{11}^{3},v_{1}^{2})
\]
\[
v_{1}^{-1}Ext_{A}^{s,t}(\mathbb{F}_{2},H_{*}(bu\wedge M))=\mathbb{F}_{2}[v_{1}^{\pm1}]
\]
In other words, the conjecture reads that $v_{1}^{-1}E_{2}(M;H)$
consists of $|H(d)|$ copies of $\mathbb{F}_{2}[v_{1}^{\pm4}]\otimes\mathbb{F}_{2}(h_{11},v_{1})/(h_{11}^{3},v_{1}^{2})$
and $|B(d)|$ copies of $\mathbb{F}_{2}[v_{1}^{\pm1}]$. To clarify,
by $|H(d)|$ we mean the number of basis elements of any given degree
in $H(d)$ and even though $H(d)$ is infinite, it is of finite type
and so for every basis element $a\in H(d)$ the copy is suspended
by the degree of $a$. The same holds for $B(d)$.

Recall $E_{3}=E_{3}(M;Q_{1})=\mathbb{F}_{2}[v_{1}^{\pm1}]\otimes\mathbb{F}_{2}[h_{11},h_{21},h_{31},...]$
with proposed differentials $d_{3}(v_{1}^{2})=h_{11}^{3}$ and $d_{3}(h_{n1})=v_{1}^{-2}h_{11}^{3}h_{n1}+v_{1}^{-2}h_{11}h_{21}h_{n-1,1}^{2}$
for $n>2$. We will express $E_{4}$ in such a way that it takes the
form Mahowald suggested. Rewrite $E_{3}=\mathbb{F}_{2}[v_{1}^{\pm1},h_{11}]\otimes\mathbb{F}_{2}[x_{1},x_{2}...]$
where $x_{n}=v_{1}h_{n+1,1}$ and introduce a grading on $E_{3}$
so that $|v_{1}^{i}|=\begin{cases}
\begin{array}{cc}
0 & \text{if }i\equiv0,1(4)\\
2 & \text{if }i\equiv2,3(4)
\end{array}\end{cases}$, $|h_{11}|=1$ and $|x_{n}|=0$. Extend this grading to monomials
in the obvious fashion. Then $E_{3}=\oplus_{n\ge0}E_{3,n}.$ The reason
we are interested in this grading is that now $d_{3}$ increases it
by $1$. But then $E_{4}$ is just the homology of the graded chain
complex i.e. $E_{4}=\oplus_{n\ge0}\ker(d_{3}^{n})/\text{im}(d_{3}^{n-1})$.

\[\begin{tikzcd} 
	0 \arrow[r,"d_3^{-1}"] & E_{3,0} \arrow[r,"d_3^0"] & E_{3,1} \arrow[r,"d_3^1"] & E_{3,2} \arrow[r,"d_3^2"] & \cdots
\end{tikzcd}\]

We claim that 

$(1)$ $\ker(d_{3}^{0})/\text{im}(d_{3}^{-1})=\ker(d_{3}^{0})=Z(d)\otimes\mathbb{F}_{2}[v_{1}^{\pm4}]\otimes\mathbb{F}_{2}[v_{1}]/(v_{1}^{2})$

$(2)$ $\ker(d_{3}^{1})/\text{im}(d_{3}^{0})=H(d)\otimes\mathbb{F}_{2}[v_{1}^{\pm4}]\otimes\mathbb{F}_{2}[v_{1}]/(v_{1}^{2})\otimes\{h_{11}\}$

$(3)$ $\begin{array}{c}
\Big(\ker(d_{3}^{2})/\text{im}(d_{3}^{1})\Big)/H(d)\otimes\mathbb{F}_{2}[v_{1}^{\pm4}]\otimes\mathbb{F}_{2}[v_{1}]/(v_{1}^{2})\otimes\{h_{11}^{2}\}\cong\\
\cong B(d)\otimes\mathbb{F}_{2}[v_{1}^{\pm4}]\otimes\mathbb{F}_{2}[v_{1}]/(v_{1}^{2})\otimes\{v_{1}^{2}\}
\end{array}$

$(4)$ $\ker(d_{3}^{n})/\text{im}(d_{3}^{n-1})=0$ for $n\ge3$

~

Given the proof of $(1)-(4)$ is not particularly insightful, we leave
it for the end of this section. We are left with the task of identifying
the expressions above with Mahowald's formulation. The key here is
to observe that given $(2)$ and $(3)$ we would need to identify
$Z(d)\otimes\mathbb{F}_{2}[v_{1}^{\pm4}]\otimes\mathbb{F}_{2}[v_{1}]/(v_{1}^{2})$
in $(1)$ with $(H(d)\oplus B(d))\otimes\mathbb{F}_{2}[v_{1}^{\pm4}]\otimes\mathbb{F}_{2}[v_{1}]/(v_{1}^{2})$.
Then from $(1),(2),(3)$ we would get the $|H(d)|$ copies of $\mathbb{F}_{2}[v_{1}^{\pm4}]\otimes\mathbb{F}_{2}(h_{11},v_{1})/(h_{11}^{3},v_{1}^{2})$.
What is left over is $B(d)\otimes\mathbb{F}_{2}[v_{1}^{\pm4}]\otimes\mathbb{F}_{2}[v_{1}]/(v_{1}^{2})$
from $(1)$ and $B(d)\otimes\mathbb{F}_{2}[v_{1}^{\pm4}]\otimes\mathbb{F}_{2}[v_{1}]/(v_{1}^{2})\otimes\{v_{1}^{2}\}$
from $(3)$, which combine to produce $|B(d)|$ copies of $\mathbb{F}_{2}[v_{1}^{\pm1}]$.
Thus each of $(1),(2)$ and $(3)$ corresponds to a third of the ``lightning
flash'' sequence, while the remainder of $(1)$ and $(3)$ each represent
half of the $v_{1}$-line.

Below we can see exactly how the elements of $H(d)$ and $B(d)$ correspond
to lightning flashes and $v_{1}$-lines in $E_{2}(M;H)$. The first
few elements of $H(d)$ appearing are $1,x_{1},x_{1}^{2},x_{2}^{2}$
and $x_{1}^{2}x_{3}+x_{2}^{3}$ and we can see the lightning falshes
for each one. Similarly, the first few elements of $B(d)$ appearing
are $x_{1}^{3}$ through $x_{1}^{9}$ and $x_{1}x_{2}^{2}$ each corresponding
to a copy of $\mathbb{F}_{2}[v_{1}^{\pm1}]$. The colors used have
no underlying meaning outside of grouping together the different elements
in $E_{2}(M;H)$ and relating each group to its representing element
of $H(d)$ or $B(d)$.\medskip{}

~

\includegraphics[scale=0.8]{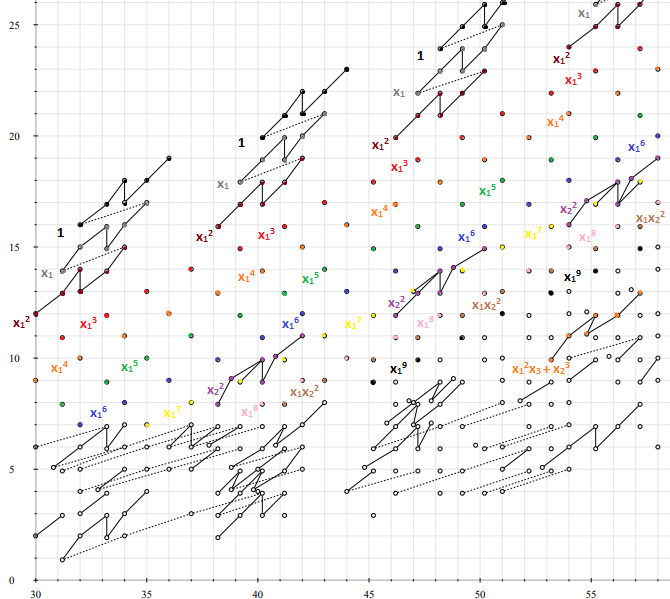}

~

We are left to prove $(1)-(4)$. It is an immediate check to verify
they follow from $(i)$ and $(ii)$ below, which is what we set out
to show.

\[
(i)\,\,\,\begin{array}{cc}
\ker(d_{3}^{n})=Z(d)\otimes\mathbb{F}_{2}[v_{1}^{\pm4}]\otimes\mathbb{F}_{2}[v_{1}]/(v_{1}^{2})\otimes\{h_{11}^{n}\} & \text{if }n=0,1\\
\ker(d_{3}^{n})/Z(d)\otimes\mathbb{F}_{2}[v_{1}^{\pm4}]\otimes\mathbb{F}_{2}[v_{1}]/(v_{1}^{2})\otimes\{h_{11}^{n}\}\cong\\
\cong B(d)\otimes\mathbb{F}_{2}[v_{1}^{\pm4}]\otimes\mathbb{F}_{2}[v_{1}]/(v_{1}^{2})\otimes\{v_{1}^{2}\}\otimes\{h_{11}^{n-2}\} & \text{if }n\ge2
\end{array}
\]

\[
(ii)\,\,\,\text{im}(d_{3}^{n})=\begin{cases}
\begin{array}{cc}
B(d)\otimes\mathbb{F}_{2}[v_{1}^{\pm4}]\otimes\mathbb{F}_{2}[v_{1}]/(v_{1}^{2})\otimes\{h_{11}^{n+1}\} & \text{if }n=0,1\\
\text{ker}(d_{3}^{n+1}) & \text{if }n\ge2
\end{array}\end{cases}
\]

Note that that $E_{3}^{0}=P\otimes\mathbb{F}_{2}[v_{1}^{\pm4}]\otimes\mathbb{F}_{2}[v_{1}]/(v_{1}^{2})$,
$E_{3}^{1}=P\otimes\mathbb{F}_{2}[v_{1}^{\pm4}]\otimes\mathbb{F}_{2}[v_{1}]/(v_{1}^{2})\otimes\{h_{11}\}$
and $d_{3}^{0}(y)=d(y)v_{1}^{-4}h_{11}$ for every $y\in P\subset E_{3}^{0}$.
Hence $\text{ker}(d_{3}^{0})$ and $\text{im}(d_{3}^{0})$ take the
desired form and the same argument holds for $\text{ker}(d_{3}^{1})$
and $\text{im}(d_{3}^{1})$. We proceed to calculate $\ker(d_{3}^{2})$
and the calculation of $\ker(d_{3}^{n})$ for $n>2$ is analogous.
Every element of $E_{3}^{2}$ takes the form $\sum_{i=1}^{s}v_{1}^{m_{i}}y_{i}+\sum_{j=1}^{t}v_{1}^{l_{j}}z_{j}h_{11}^{2}$
where $m_{1}<m_{2}<\cdots<m_{s}$, $m_{i}\equiv2,3(4)$, $l_{1}<l_{2}<\cdots l_{t},$
$l_{j}\equiv0,1(4)$ and $y_{i},z_{j}\in P$. We also assume $y_{i},z_{j}\neq0$.
Then 

\[
d_{3}^{2}\Bigg(\sum_{i=1}^{s}v_{1}^{m_{i}}y_{i}+\sum_{j=1}^{t}v_{1}^{l_{j}}z_{j}h_{11}^{2}\Bigg)=\sum_{i=1}^{s}\big(v_{1}^{m_{i}-2}y_{i}h_{11}^{3}+v_{1}^{m_{i}-4}d(y_{i})h_{11}\big)+\sum_{j=1}^{t}v_{1}^{l_{j}-4}d(z_{j})h_{11}^{3}
\]

Setting this equal to $0$ we observe two cases. First if $s=0$ then
$d(z_{j})=0$ for all $j$ and we get the same component as in $\ker(d_{3}^{0})$,
namely $Z(d)\otimes\mathbb{F}_{2}[v_{1}^{\pm4}]\otimes\mathbb{F}_{2}[v_{1}]/(v_{1}^{2})\otimes\{h_{11}^{2}\}\subset\ker(d_{3}^{2})$.
If $s>0$ then we obtain $d(y_{i})=0$ for all $i$ and we are left
with 
\[
\sum_{i=1}^{s}v_{1}^{m_{i}-2}y_{i}+\sum_{j=1}^{t}v_{1}^{l_{j}-4}d(z_{j})=0
\]

\noindent which given the degrees of $v_{1}$ can only happen if $s=t$,
$m_{i}-2=l_{i}-4$ and $y_{i}=d(z_{i})$. Note $y_{i}=d(z_{i})$ already
implies $d(y_{i})=0$. Furthermore, for every $y_{i}\in B(d)$ we
have a unique $z_{i}\in P$ with $y_{i}=d(z_{i})$ modulo $Z(d)\otimes\mathbb{F}_{2}[v_{1}^{\pm4}]\otimes\mathbb{F}_{2}[v_{1}]/(v_{1}^{2})\otimes\{h_{11}^{2}\}\subset\ker(d_{3}^{2})$.
Hence 
\[
\ker(d_{3}^{2})/Z(d)\otimes\mathbb{F}_{2}[v_{1}^{\pm4}]\otimes\mathbb{F}_{2}[v_{1}]/(v_{1}^{2})\otimes\{h_{11}^{2}\}\cong B(d)\otimes\mathbb{F}_{2}[v_{1}^{\pm4}]\otimes\mathbb{F}_{2}[v_{1}]/(v_{1}^{2})\otimes\{v_{1}^{2}\}
\]

\noindent as desired. In fact, $\ker(d_{3}^{2})\cong P\otimes\mathbb{F}_{2}[v_{1}^{\pm4}]\otimes\mathbb{F}_{2}[v_{1}]/(v_{1}^{2})$,
but stated this way it does not relate well with Mahowald's conjecture.

Next we show $\text{im}(d_{3}^{2})=\ker(d_{3}^{3})$ and the result
for $\text{im}(d_{3}^{n})$ follows analogically. As we saw above
elements of $\ker(d_{3}^{3})$ are sums of elements of the form $v_{1}^{m}yh_{11}+v_{1}^{m-2}zh_{11}^{3}$
for $m\equiv2,3(4)$ and $y,z\in P$ such that $d(z)=y$. But then
$d_{3}^{2}(v_{1}^{m}z)=v_{1}^{m}yh_{11}+v_{1}^{m-2}zh_{11}^{3}$ and
so $\ker(d_{3}^{3})\subset\text{im}(d_{3}^{2})$ and since the reverse
inclusion holds as well the two must coincide. This completes the
proof of $(i)$ and $(ii)$ and thus we have successfully identified
Mahowald's and Palmieri's formulations of the problem.


\begin{thebibliography}{1}
\bibitem{key-1}M. Mahowald, $\emph{The order of the image of the J-homomorphism}$,
Bull. Amer. Math. Soc. 76 (1970), 1310-1313

\bibitem{key-2}J. Palmieri, $\emph{Stable homotopy over the Steenrod algebra}$,
Amer. Math. Soc. (2001) Vol. 151 N. 716\end{thebibliography}
\end{document}